\documentclass{article}

\usepackage{amsmath,amsfonts,amsthm,amssymb,graphicx,tikz,array}
\usepackage[all,2cell,ps]{xy}

\theoremstyle{plain}
\newtheorem{thm}{Theorem}

\newtheorem{lem}[thm]{Lemma}
\newtheorem{prop}[thm]{Proposition}

\newtheorem{cor}[thm]{Corollary}

\theoremstyle{definition}

\newtheorem*{rem}{Remark}

\theoremstyle{remark}

\newcommand{\bbB}{\mathbb{B}}
\newcommand{\bbC}{\mathbb{C}}

\newcommand{\bbQ}{\mathbb{Q}}

\newcommand{\bbZ}{\mathbb{Z}}

\newcommand{\calA}{\mathcal{A}}

\newcommand{\calI}{\mathcal{I}}
\newcommand{\calJ}{\mathcal{J}}

\newcommand{\calM}{\mathcal{M}}
\newcommand{\calN}{\mathcal{N}}
\newcommand{\calO}{\mathcal{O}}

\newcommand{\calZ}{\mathcal{Z}}

\newcommand{\al}{\alpha}
\newcommand{\gam}{\gamma}
\newcommand{\Gam}{\Gamma}

\newcommand{\de}{\delta}

\newcommand{\lam}{\lambda}
\newcommand{\Lam}{\Lambda}
\newcommand{\sig}{\sigma}

\DeclareMathOperator{\M}{M}

\DeclareMathOperator{\PSL}{PSL}

\DeclareMathOperator{\PU}{PU}

\DeclareMathOperator{\PSU}{PSU}

\DeclareMathOperator{\Isom}{Isom}

\newcommand{\ssm}{\smallsetminus}

\newcommand{\conj}{\overline}
\newcommand{\wh}{\widehat}

\newenvironment{pf}{\begin{proof}}{\end{proof}}

\newenvironment{enum}{\begin{enumerate}}{\end{enumerate}}

\usetikzlibrary{arrows}

\title{On general type surfaces with $q=1$ and $c_2 = 3 p_g$}
\author{Matthew Stover\footnote{This material is based upon work supported by the National Science Foundation under Grant Number NSF 1361000 and Grant Number 523197 from the Simons Foundation/SFARI.
} \\ \small{Temple University}\\ \small{\textsf{mstover@temple.edu}}}
\date{\today}

\begin{document}

\maketitle

\begin{abstract}
Let $S$ be a minimal surface of general type with irregularity $q(S) = 1$. Well-known inequalities between characteristic numbers imply that
\[
3 p_g(S) \le c_2(S) \le 10 p_g(S),
\]
where $p_g(S)$ is the geometric genus and $c_2(S)$ the topological Euler characteristic. Surfaces achieving equality for the upper bound are classified, starting with work of Debarre. We study equality in the lower bound, showing that for each $n \ge 1$ there exists a surface with $q = 1$, $p_g = n$, and $c_2 = 3n$. The moduli space $\calM_n$ of such surfaces is a finite set of points, and we prove that $\#\calM_n \to \infty$ as $n \to \infty$. Equivalently, this paper studies the number of closed complex hyperbolic $2$-manifolds of first betti number $2$ as a function of volume; in particular, such a manifold exists for every possible volume.
\end{abstract}

\section{Introduction}\label{sec:intro}

One of the primary problems in the study of algebraic surfaces is to classify the smooth minimal surfaces of general type with given characteristic numbers. We refer to \cite{BCP1} for a survey on this problem. The first result of this paper is the following.

\begin{thm}\label{thm:exists}
For every $n \ge 1$, there exists a smooth minimal complex projective surface $S_n$ of general type with irregularity $q(S_n) = 1$, geometric genus $p_g(S_n) = n$, and topological Euler number $c_2(S_n) = 3 n$.
\end{thm}

We also study the moduli space $\calM_n$ of such surfaces. A surface satisfying the conditions of Theorem \ref{thm:exists} is necessarily a ball quotient, so $\calM_n$ is finite and, by Theorem \ref{thm:exists}, nonempty. Recall that the minimal smooth complex projective surfaces of general type satisfying $c_1^2 = 3 c_2$ are all of the form $\bbB^2 / \Gam$ with $\bbB^2$ the unit ball in $\bbC^2$ and $\Gam$ a torsion-free cocompact lattice in $\PU(2,1)$. Finiteness of $\calM_n$ then follows from Mostow--Siu rigidity \cite{Siu}. We also give lower and upper bounds for $\# \calM_n$ as a function of $n$.

\begin{thm}\label{thm:moduli}
For $n \ge 1$, let $\calM_n$ be the moduli space of surfaces satisfying the conditions of Theorem \ref{thm:exists}. Then $\calM_n$ is a finite nonempty set of points. Moreover, there are universal constants $C_1, C_2 > 0$ such that
\begin{equation}\label{eq:bounds}
C_1\, n < \# \calM_n < e^{C_2 n}
\end{equation}
for all $n$. In particular, $\# \calM_n \to \infty$ as $n \to \infty$.
\end{thm}

We also make a comment on the analogous problem for $q = 0$ and $c_2 = 3 p_g + 3$ at the end of the paper.
The upper bound follows from a theorem of Gelander \cite[Thm.\ 1.11]{Gelander}. We note that algebraic results of Catanese imply the upper bound $6^{(9 n + 5/9)^{15}}$ \cite[Thm.\ A]{Catanese}. See \cite{LP, LP2} for a lower bound for the number of general type surfaces with given $c_1^2$; for our surfaces $c_1^2 = 9 n$. Our method of proof for the lower bound is explicit. One can take $S_1$ to be the Cartwright--Steger surface \cite{CS} and our surfaces $S_n$ are \'etale abelian covers of $S_1$ of degree $n$. To prove that our surfaces satisfy the conditions of Theorem \ref{thm:exists}, we must show that $q(S_n) = 1$. This follows from the following result.

\begin{thm}\label{thm:alexander}
Any finite \'etale abelian cover of the Cartwright--Steger surface has irregularity $q = 1$. In other words, there is no jumping in first cohomology for abelian covers.
\end{thm}

Let $\Gam$ be the fundamental group of the Cartwright--Steger surface and identify the group ring of its abelianization $\Gam^{ab} \cong \bbZ^2$ with $\bbZ[r^{\pm 1}, s^{\pm 1}]$. To prove Theorem \ref{thm:alexander}, we will show that for our chosen generating set for $\Gam^{ab}$ the Alexander stratification in the sense of Hironaka \cite{Hironaka} (see \S \ref{ssec:Alexander}) is:
\[
\wh{\Gam} \supset V_0(\Gam) = \{r - s\} \supset V_1(\Gam) = \{\,\wh{1}\,\} \supset V_2(\Gam) = \emptyset
\]
Here $\wh{\Gam}$ denotes the character group of $\Gam$ and $\wh{1}$ the trivial character. Then $r - s$ represents the $\bbC^* \subset \wh{\Gam}$ uniquely determined by sending the two generators for $\Gam^{ab}$ to the same element of $\bbC^*$. That there are no cohomology jumps in finite \'etale abelian covers follows directly from the fact that the only finite character in any $V_i(\Gam)$ is $\wh{1}$.

However, to show that the number of choices of $S_n$ grows linearly with $n$, we must count the number of \emph{nonisomorphic} abelian covers of $S_1$ with given degree $n$. In particular, one must take care of the fact that inequivalent covers of $S_1$ may give isomorphic surfaces. Indeed, two nonconjugate finite index subgroups of $\Gam$ may well be conjugate in $\Isom(\bbB^2)$ and therefore determine biholomorphic ball quotients. Mostow--Siu rigidity \cite{Siu} implies that two closed ball quotient manifolds are biholomorphic (in fact, homeomorphic) if and only if their fundamental groups are conjugate in $\Isom(\bbB^2)$, so proving Theorem \ref{thm:moduli} is a counting problem for conjugacy classes of lattices in the Lie group $\Isom(\bbB^2)$. Restating our results in this language, we have the following.

\begin{cor}\label{cor:PUversion}
Let $\calM_n$ be the set of isomorphism classes of torsion-free cocompact lattices $\Gam$ in $\PU(2,1)$ with first betti number $b_1(\Gam) = 2$ and Euler characteristic $e(\Gam) = 3n$. Then $\calM_n$ is nonempty for all $n$ and there are universal constants $C_1, C_2 > 0$ such that $\# \calM_n$ satisfies \eqref{eq:bounds}. In particular, for every possible volume of a closed complex hyperbolic $2$-manifold, there is a manifold of that volume and first betti number exactly $2$.
\end{cor}

The claim about volume is immediate from Hirzebruch proportionality \cite{HirzebruchProportionality} and Chern--Gauss--Bonnet. In fact, for all $n \ge 1$, one can find $\Gam_n \in \calM_n$ such that $\{\Gam_{k n}\}_{k \in \bbZ}$ is a nested family of lattices with $\Gam_{k n} \in \calM_{k n}$ (i.e., $\Gam_{(k+1) n} \subset \Gam_{k n}$ for all $k$). We close with one final immediate consequence of our work; see \cite{FV} for some interest in problems of this kind.

\begin{cor}\label{cor:Alb}
There are infinitely many $2$-dimensional smooth ball quotients of Albanese dimension $1$. In fact, there is one achieving every possible volume. One may take the infinite collection to lie in a tower of finite \'etale abelian covers.
\end{cor}

As a final remark in this direction, we learned after completing this paper that very recent work of Vidussi proves Theorem \ref{thm:alexander} for certain cyclic coverings \cite{Vidussi}. More precisely, he proves that cyclic covers of the Cartwright--Steger surface of degree $d = \lam e + 1$ have irregularity one, where $\lam \in \bbZ$ and $e$ is the least common multiple of the orders of the elements in the Green--Lazarsfeld set $W_1$. He then uses ramified double covers of these cyclic coverings to produce smooth minimal surfaces of general type with Chern slopes dense in the interval $[8,9]$. In particular, \cite{Vidussi} gives another very interesting application of the surfaces studied in this paper.

\medskip

\noindent\textbf{Acknowledgments.} Many thanks are due to Fabrizio Catanese for suggesting that I consider the jumping loci for the Cartwright--Steger surface, as well as for some comments on the previous literature. I also thank the referee for helpful suggestions.

\section{Preliminaries}\label{sec:prelims}

\subsection{The basic inequality}

Let $S$ be a smooth minimal complex projective surface of general type with irregularity $q = 1$ and geometric genus $p_g$. One immediately obtains that the holomorphic Euler characteristic is
\begin{equation}\label{eq:chi}
\chi = p_g.
\end{equation}
Let $c_1^2 = K_S^2$ be the self-intersection of the canonical divisor and $c_2$ be the topological Euler characteristic. We also have Noether's formula
\begin{equation}\label{eq:noether}
12 \chi = c_1^2 + c_2
\end{equation}
along with the Hodge decomposition
\begin{equation}\label{eq:hodge}
c_2 = 2 - 4 q + 2 p_g + h^{1,1}.
\end{equation}
We then have the following pair of inequalities:

\begin{lem}\label{lem:MainIneq}
Suppose that $S$ is a smooth minimal complex projective surface of general type with irregularity $q = 1$. Then the geometric genus $p_g$ and topological Euler characteristic $c_2$ satisfy
\begin{equation}\label{eq:MainIneq}
3 p_g \le c_2 \le 10 p_g.
\end{equation}
Moreover, $S$ achieves equality for the lower bound if and only if $S$ is a ball quotient.
\end{lem}

\begin{pf}
We will show that the first inequality is equivalent to the Bogomolov--Miyaoka--Yau inequality $c_1^2 \le 3 c_2$. Combining this with \eqref{eq:chi} and \eqref{eq:noether} we have
\[
c_1^2 = 12 p_g - c_2 \le 3 c_2,
\]
and the lower bound is immediate. Since $c_1^2 = 3 c_2$ if and only if $S$ is a ball quotient, the last assertion of the lemma also follows. On the other hand, Debarre proved that $q > 0$ implies that $c_1^2 \ge 2 p_g$ \cite{Debarre}. Therefore \eqref{eq:chi} and \eqref{eq:noether} now give
\[
12 p_g - c_2 \ge 2 p_g,
\]
which proves the upper bound.
\end{pf}

\begin{rem}
As mentioned in the introduction, the classification of surfaces achieving equality in the upper bound of \eqref{eq:MainIneq} was very recently completed. See \cite{CMLP}.
\end{rem}

\subsection{Alexander stratifications and cohomology jumps}\label{ssec:Alexander}

See \cite{Hironaka} for an excellent treatment of the material in this section. Let $\Gam = \langle g_i\ |\ R_j \rangle$ be a finitely presented group with abelianization $\Gam^{ab}$ and $\al : \Gam \to \Gam^{ab}$ be the abelianization. One then defines the \emph{Fox derivative} $D_i = \partial / \partial g_i$, which maps $\Gam$ to its group ring $\bbZ[\Gam^{ab}]$ by the rules:
\begin{align}
D_i(g_j) &= \de_{i,j} \label{eq:FoxRule1} \\
D_i(g h) &= D_i(g) + \al(g) D_i(h) \label{eq:FoxRule2}
\end{align}
If $\Gam$ has $n$ generators and $m$ relations, we then obtain the \emph{Alexander matrix}, which is the $n \times m$ matrix $\calA(\Gam) = (D_i(R_j))$ with coefficients in $\bbZ[\Gam^{ab}]$. One has the following algorithm to compute the Alexander matrix.

\begin{lem}\label{lem:AlexanderAlgorithm}
Let $g_i$ be a generator of the group $\Gam$ with presentation $\langle g_i\ |\ R_j \rangle$ and $\al$ the abelianization homomorphism to the group ring $\bbZ[\Gam^{ab}]$. If
\[
R_j = g_{i_1}^{\ell_1} \cdots g_{i_t}^{\ell_t}
\]
is a relation, then the following algorithm computes the Fox derivative $\partial R_j / \partial g_i$:
\begin{enum}

\item Remove all generators in $R_j$ to the right of the last appearance of $g_i$.

\item For $i_k \neq i$, replace $g_{i_k}^{\ell_k}$ with $\al(g_{i_k})^{\ell_k}$.

\item Replace any appearance of $g_i^\ell$ with $(D_i(g_i^\ell) + \al(g_i)^\ell)$.

\item Replace $D_i(g_i^\ell)$ with:
\begin{enumerate}
\item $\displaystyle{\sum_{j = 0}^{\ell - 1} \al(g_i)^j \quad (\ell \ge 1)}$
\item $\displaystyle{\sum_{j = -1}^{\ell} \al(g_i)^j \quad (\ell \le -1)}$
\end{enumerate}

\item Simplify the polynomial.

\end{enum}
\end{lem}

\begin{pf}
To prove that \emph{1.}\ is valid, we want to show that if we divide the product decomposition of $R_j$ into $\gam_1 \gam_2$, where $g_i$ does not appear in $\gam_2$, then $D_i(R_j) = D_i(\gam_1)$. However, induction on \eqref{eq:FoxRule1} and \eqref{eq:FoxRule2} gives $D_i(\gam_2) = 0$. Then we have
\[
D_i(R_j) = D_i(\gam_1) + \al(\gam_1) D_i(\gam_2) = D_i(\gam_1),
\]
so \emph{1.}\ holds.

For \emph{2.}, suppose that we have $R_j = \gam_1 g_{i_k}^{\ell_{i_k}} \gam_2$ with $i_k \neq i$. Then $D_i(R_j)$ equals
\begin{align*}
&D_i(\gam_1) + \al(\gam_1)D_i(g_{i_k}^{\ell_{i_k}} \gam_2) \\
= &D_i(\gam_1) + \al(\gam_1)(D_i(g_{i_k}^{\ell_{i_k}}) + \al(g_{i_k})^{\ell_k} D_i(\gam_2)) \\
= &D_i(\gam_1) + \al(\gam_1) \al(g_{i_k})^{\ell_{i_k}} D_i(\gam_2) \\
= &(D_i(g_{i_1}^{\ell_1}) + \al(g_{i_1}^{\ell_1})(D_i(g_{i_2}^{\ell_2}) + \al(g_{i_2}^{\ell_2})( \cdots \al(g_{i_{k-1}})^{\ell_{i_{k-1}}} \al(g_{i_k})^{\ell_{i_k}} D_i(\gam_2)) \cdots )).
\end{align*}
Similarly, we have
\begin{align*}
D_i(\gam_1 \gam_2) &= D_i(\gam_1) + \al(\gam_1)D_i\gam_2) \\
&= D_i(\gam_1) + \al(\gam_1) D_i(\gam_2) \\
&= (D_i(g_{i_1}^{\ell_1}) + \al(g_{i_1}^{\ell_1})(D_i(g_{i_2}^{\ell_2}) + \al(g_{i_2}^{\ell_2})( \cdots \al(g_{i_{k-1}})^{\ell_{i_{k-1}}} D_i(\gam_2)) \cdots )),
\end{align*}
so one obtains $D_i(R_i)$ from $D_i(\gam_1 \gam_2)$ by inserting $\al(g_{i_k})^{\ell_{i_k}}$ in $D_i(\gam_1 \gam_2)$ between $\al(g_{i_{k-1}})^{\ell_{i_{k-1}}}$ and $D_i(\gam_2)$, which is precisely what \emph{2.}\ does. The argument that \emph{3.}\ is valid is exactly the same.

Finally, \emph{4.}\ is an easy induction on \eqref{eq:FoxRule2} and \emph{5.}\ is just bookkeeping.
\end{pf}

Let $\wh{\Gam}$ be the character group of $\Gam$, and $\wh{1}$ will denote the trivial character. Considering $\bbZ[\Gam^{ab}]$ as a ring of Laurent polynomials, any $\rho \in \wh{\Gam}$ defines an `evaluation map' $\bbZ[\Gam^{ab}] \to \bbC$ in a canonical way \cite[\S 2.1]{Hironaka}. In particular, we can consider the $n \times m$ matrix
\[
\calA_\rho(\Gam) \in \M_{n \times m}(\bbC)
\]
determined by evaluating the Alexander matrix $\calA(\Gam)$ at $\rho$ and then define
\begin{equation}\label{eq:JumpingDef}
V_i(\Gam) = \left\{ \rho \in \wh{\Gam}\ |\ \mathrm{rank}(\calA_\rho(\Gam)) < n - i \right\}
\end{equation}
for $0 \le i < n$. We call $\{V_i(\Gam)\}$ the \emph{Alexander stratification} of $\Gam$.

Now, let $G$ be a finite abelian group such that there exists a surjective homomorphism $\al : \Gam \to G$. We then obtain a natural embedding of character groups
\[
\wh{\al} : \wh{G} \to \wh{\Gam}.
\]
Considering $f \in \bbZ[\Gam^{ab}]$ as a Laurent polynomial in variables $x_1, \dots, x_r$ (so $\Gam^{ab} \cong \bbZ^r$ modulo torsion), notice that the value of $f$ on the character $\al$ of $G$ is given by evaluating the Laurent polynomial at roots of unity associated with the cyclic subgroups of $G$ generated by the images of the fixed generators of $\Gam^{ab}$.

For a finitely generated group $\Lam$, let $b_1(\Lam)$ denote the first betti number of $\Lam$, i.e., the rank of $\Lam^{ab} \otimes_\bbZ \bbQ$. We then have the following.

\begin{prop}[Prop.\ 2.5.6 \cite{Hironaka}]\label{prop:Jumping}
Let $\Gam$ be a finitely presented group on $n$ generators with Alexander matrix $\calA(\Gam)$. If $\al : \Gam \to G$ is a surjective homomorphism onto a finite abelian group $G$, let $\Gam_\al$ be the kernel of $\al$. Then
\begin{equation}\label{eq:Jumping}
b_1(\Gam_\al) = b_1(\Gam) + \sum_{i = 1}^{n - 1} \left|V_i(\Gam) \cap \wh{\al}(\wh{G} \ssm \{\wh{1}\})\right|.
\end{equation}
In particular, if each $V_i(\Gam)$ contains no finite characters other than possibly the trivial character, then $b_1(\Gam_\al) = b_1(\Gam)$.
\end{prop}

\section{Proofs of Theorems \ref{thm:exists}-\ref{thm:alexander}}\label{sec:proofs1}

We rely heavily on the notation from \S \ref{ssec:Alexander}. We now jump directly into proving Theorem \ref{thm:alexander}.

\begin{pf}[Proof of Theorem \ref{thm:alexander}]
The fundamental group $\Gam$ of the Cartwright--Steger surface has generators $x,y,z$ and relations:
\begin{align*}
R_1 =\,\, & y^{-1} z^{-2} x^{-3} z^{-1} x^{-1} z y^{-2} z^{-1} y x^{-3} \\
R_2 =\,\, & y^3 x^3 z y^{-1} z^{-1} x^{-3} y^{-2} z x z^{-1} x^{-1} \\
R_3 =\,\, & y z^{-1} x^{-1} z^{-1} x^{-3} y^{-3} z^{-1} x^{-3} z^{-1} y^{-1} z y \\
R_4 =\,\, & y z x y^3 x^3 z y^{-1} x^3 z^2 y^{-1} z^{-1} \\
R_5 =\,\, & z^{-1} x^{-3} y^{-3} z x^{-1} z^{-1} y^{-1} z x y^3 x^3 y \\
R_6 =\,\, & z^{-1} x^{-3} y^{-2} z^{-2} x^{-3} y^{-1} z x z^{-1} y^{-1} z y^2 z^{-1} x^{-2} \\
R_7 =\,\, & z x^{-1} z^{-1} y^2 x^3 z^2 x z^{-1} y^{-1} z y^2 z^{-2} x^{-3} y^{-3} \\
R_8 =\,\, & y x^2 z y^{-2} z^{-1} y x^{-3} z^{-2} x^{-3} y^{-3} z y z^{-2} x^{-3} \\
R_9 =\,\, & y^{-2} x y^3 z^{-1} y^{-1} z y^2 z^{-1} x^{-1} y^{-1} z x^3 y^{-1} z y^{-1} z y z^{-2} x^{-3} \\
R_{10} =\,\, & z^{-1} x^3 z^2 y^{-1} z^{-1} y^4 x^3 z y^{-1} z y^3 x^3 z y^{-1} z x z^{-1} x^{-3} y^{-2} \\
R_{11} =\,\, & z^{-1} x^{-3} y^{-3} z^{-1} y z y^{-2} z^{-1} y^4 x^3 z y^{-1} z x z^{-1} x^{-3} z^{-1} \\
&x^{-2} z y^{-1} z^{-2} x^{-3} \\
R_{12} =\,\, & y^{-1} z y z^{-2} x^{-3} y^3 x^3 z x^2 z x^3 z x^{-1} z^{-1} y z^{-1} x^{-3} y^{-3} z^{-1} y \\
&z y^{-2} z^{-1} y z x^{-1} z^{-1} y^{-1} z
\end{align*}
Let $\al : \Gam \to \bbZ^2$ be the abelianization, $a = \al(x)$, $b = \al(y)$, and $c = \al(z)$. Then each relation $R_i$ becomes either trivial or equivalent to the relation
\[
7 a + 2 b + 3 c = 0.
\]
Considering $\bbZ[r^{\pm 1}, s^{\pm 1}]$ as the group ring of $\bbZ^2$, we can identify $a$ with $r s^{-1}$, $b$ with $r s^2$, and $c$ with $r^{-3} s$. The Alexander matrix $\calA(\Gam)$ is then the $3 \times 12$ matrix determined by the entries in Tables \ref{tb:ddx} - \ref{tb:ddz}.

One can then check directly with a computer algebra program that all $3 \times 3$ minors of $\calA(\Gam)$ have determinant zero if and only if $r = s$. Furthermore, $\calA(\Gam)$ always has rank $2$ when $r = s$ except for precisely the case when $r = s = 1$. Thus the Alexander stratification for $\Gam$ is:
\[
\wh{\Gam} \supset V_0(\Gam) = \{r - s\} \supset V_1(\Gam) = \{\,\wh{1}\,\} \supset V_2(\Gam) = \emptyset
\]
Proposition \ref{prop:Jumping} immediately implies the conclusion of the theorem. Indeed,
\[
V_i \ssm \wh{\al}(\wh{G} \ssm \{\wh{1}\}) = \emptyset
\]
for $i = 1, 2$, so the sum on the right hand side of \eqref{eq:Jumping} is zero, hence $b_1(\Gam_\al) = b_1(\Gam) = 2$ for all $\al$.
\end{pf}

\pagebreak

\begin{center}
\vspace*{\stretch{1}}

\begin{table}[ht]
\begin{centering}
\begin{tabular}{| >{$}>{\tiny}c<{$} | >{$}>{\tiny}c<{$} |}
\hline
 & \\[-0.6em]
\partial R_1 / \partial x & -\frac{r^4 s^2+r^4+r^3 s+r^2 s^2+r^2 s+r s^2+s^3}{s^3} \\
 & \\[-0.6em]
\hline
 & \\[-0.6em]
\partial R_2 / \partial x & \frac{r^7 s^5+r^6 s^6-r^6 s^3+r^5 s^7-r^5 s^4-r^4 s^5-r^3+s}{r^2 s} \\
 & \\[-0.6em]
\hline
 & \\[-0.6em]
\partial R_3 / \partial x & -\frac{r \left(r^4 s^5+r^3 s^6+r^2 s^7+r^2 s^5+r^2+r s+s^2\right)}{s^3} \\
 & \\[-0.6em]
\hline
 & \\[-0.6em]
\partial R_4 / \partial x & \frac{s^2 \left(r^6 s^4+r^5 s^5+r^5+r^4 s^6+r^4 s+r^3 s^2+s\right)}{r^2} \\
 & \\[-0.6em]
\hline
 & \\[-0.6em]
\partial R_5 / \partial x & -\frac{r^9 s^5+r^8 s^6+r^7 s^7-r^6 s^3-r^5 s^4-r^4 s^5+r s^2-1}{r^7 s^5} \\
 & \\[-0.6em]
\hline
 & \\[-0.6em]
\partial R_6 / \partial x & -\frac{r^5+2 r^4 s+r^3 s^4+2 r^3 s^2+r^2 s^5+r s^6-1}{r s^4} \\
 & \\[-0.6em]
\hline
 & \\[-0.6em]
\partial R_7 / \partial x & -\frac{s^2 \left(r^9 s^2+r^8 s^3+r^7 s^4-r^7 s-r^6 s^2-r^5 s^3-r^2 s^2+1\right)}{r^4} \\
 & \\[-0.6em]
\hline
 & \\[-0.6em]
\partial R_8 / \partial x & -\frac{r^6 s^3+r^5 s^4+r^5+r^4 s^5+r^4 s-r^3 s^5+r^3 s^2-r^2 s^6+r^2 s^3+r s^4+s^5}{r s^4} \\
 & \\[-0.6em]
\hline
 & \\[-0.6em]
\partial R_9 / \partial x & -\frac{r^7 s^7-r^5 s^4-r^4 s^5+r^4 s^2-r^3 s^6+r^3 s^3+r^2 s^4-1}{r^2 s^4} \\
 & \\[-0.6em]
\hline
 & \\[-0.6em]
\partial R_{10} / \partial x & \frac{r \left(r^7 s^4+r^6 s^7+r^6 s^5+r^5 s^8+r^5 s^6+r^4 s^9+r^4-r^3 s^5+r^3 s-r^2 s^6+r^2 s^2-r s^7+s^6\right)}{s^3} \\
 & \\[-0.6em]
\hline
 & \\[-0.6em]
\partial R_{11} / \partial x & \frac{r^6 s^2+r^5 s^3+r^4 s^4-r^4 s-r^4-r^3 s^2-r^3 s-r^3-r^2 s^4-r^2 s^2-r^2 s-r s^5-r s^2-s^6+s}{s^4} \\
 & \\[-0.6em]
\hline
 & \\[-0.6em]
\partial R_{12} / \partial x & \scriptstyle{-r^8 s^4-r^7 s^5-r^6 s^6+r^6 s^5+r^5 s^6+r^5 s^4+r^5 s^2+r^4 s^7+r^4 s^5+r^4 s^3+r^3 s^4-r^3-r^2 s^3-r^2-r s-s^2} \\[-0.6em]
 & \\
\hline
\end{tabular}
\caption{$\partial R_i / \partial x$}\label{tb:ddx}
\end{centering}
\end{table}

\vspace*{\stretch{1}}
\end{center}

\clearpage

\pagebreak

\begin{center}
\vspace*{\stretch{1}}

\begin{table}[ht]
\begin{centering}
\begin{tabular}{| >{$}>{\small}c<{$} | >{$}>{\small}c<{$} |}
\hline
 & \\[-0.6em]
\partial R_1 / \partial y & \frac{r^3-r s^3-s^3-s}{r s^5} \\
 & \\[-0.6em]
\hline
 & \\[-0.6em]
\partial R_2 / \partial y & \scriptstyle{r^2 (s-1) s^2 (s+1)} \\
 & \\[-0.6em]
\hline
 & \\[-0.6em]
\partial R_3 / \partial y & -\frac{r^3 s^5+r^3+r^2 s^3-r s^3+r s-s}{r s^3} \\
 & \\[-0.6em]
\hline
 & \\[-0.6em]
\partial R_4 / \partial y & \frac{r^4 s^6-r^4 s^4+r^3 s^4+r^3+r^2 s^2-s}{r^3} \\
 & \\[-0.6em]
\hline
 & \\[-0.6em]
\partial R_5 / \partial y & -\frac{r^6 s^5-r^6 s^3+r^5 s^3+r^4 s-r^2 s^4+r^2-r s^2-1}{r^7 s^5} \\
 & \\[-0.6em]
\hline
 & \\[-0.6em]
\partial R_6 / \partial y & -\frac{r^3 s^3+r^3+r^2 s^6+r s^4-r s^3-s}{r^3 s^6} \\
 & \\[-0.6em]
\hline
 & \\[-0.6em]
\partial R_7 / \partial y & -\frac{r^4 s^4+r^3 s^2+r^3-r^2 s^3+r^2-r s^3-r s-s}{r^2} \\
 & \\[-0.6em]
\hline
 & \\[-0.6em]
\partial R_8 / \partial y & -\frac{r^5 s^4-r^5-r^4 s^4+r^4 s^2+r^3 s^3+r^3+r^2 s-s}{r^4 s^4} \\
 & \\[-0.6em]
\hline
 & \\[-0.6em]
\partial R_9 / \partial y & -\frac{r^8 s^6+r^8 s^3+r^7 s^2-r^6 s^6-2 r^5 s^4-r^4 s^2+r^3 s^3+r^3 s-r^3+r^2 s-s^2}{r^4 s^5} \\
 & \\[-0.6em]
\hline
 & \\[-0.6em]
\partial R_{10} / \partial y & \frac{r^6 s^6-r^6 s^4+r^5 s^9-r^5 s^7+r^5 s^4+r^4 s^7+r^4 s^2+r^3 s^5+r^3-r^2 s^7-r s^5-s}{r s^5} \\
 & \\[-0.6em]
\hline
 & \\[-0.6em]
\partial R_{11} / \partial y & -\frac{-r^6 s^6+r^6 s^4-r^5 s^4-2 r^4 s^2+r^3 s^7-r^3+r^2 s^5+r s^6+2 r s^3+s}{r^4 s^7} \\
 & \\[-0.6em]
\hline
 & \\[-0.6em]
\partial R_{12} / \partial y & \frac{r^9 s^2-r^8 s^7+r^8 s^5+r^8-r^7 s^5-r^7 s+r^6 s^8-2 r^6 s^3+r^5 s^6-r^5 s+r^4 s^4-r^3+s}{r^4 s^2} \\[-0.6em]
 & \\
\hline
\end{tabular}
\caption{$\partial R_i / \partial y$}\label{tb:ddy}
\end{centering}
\end{table}

\begin{table}[ht]
\begin{centering}
\begin{tabular}{| >{$}>{\small}c<{$} | >{$}>{\small}c<{$} |}
\hline
 & \\[-0.6em]
\partial R_1 / \partial z & -\frac{r^2 \left(r^3 s^3+r^3 s-r^2 s^4+s^2+1\right)}{s^5} \\
 & \\[-0.6em]
\hline
 & \\[-0.6em]
\partial R_2 / \partial z & \frac{r^6 s^4-r^5 s^2-r+s}{s} \\
 & \\[-0.6em]
\hline
 & \\[-0.6em]
\partial R_3 / \partial z & -\frac{r^2 \left(r^4 s^4+r^2 s^4+r s^2+r-1\right)}{s^3} \\
 & \\[-0.6em]
\hline
 & \\[-0.6em]
\partial R_4 / \partial z & \scriptstyle{r^5 s^5+r^4 s+2 r s^2-1} \\
 & \\[-0.6em]
\hline
 & \\[-0.6em]
\partial R_5 / \partial z & -\frac{r^8 s^4-r^2 s+r s^2-1}{r^5 s^5} \\
 & \\[-0.6em]
\hline
 & \\[-0.6em]
\partial R_6 / \partial z & -\frac{r^4 s^2+r^3 s^5+r^2 s^4+r s^3+r s^2-s^3-1}{s^6} \\
 & \\[-0.6em]
\hline
 & \\[-0.6em]
\partial R_7 / \partial z & -\frac{r^7 s^3-r^5 s^2+r^4 s^4+r^3 s^2-r^2 s^3-r^2-r+s}{r} \\
 & \\[-0.6em]
\hline
 & \\[-0.6em]
\partial R_8 / \partial z & -\frac{r^6 s^3-r^4 s^4+r^4 s+r^3 s^4+r^2+r s^2-1}{r s^4} \\
 & \\[-0.6em]
\hline
 & \\[-0.6em]
\partial R_9 / \partial z & -\frac{r^7 s^6+r^6 s^4-r^5 s^5-r^5 s^2+r s^2-1}{r s^4} \\
 & \\[-0.6em]
\hline
 & \\[-0.6em]
\partial R_{10} / \partial z & \frac{r^2 \left(r^7 s^5+r^6 s^8+r^4 s-r^3 s^6+r^3 s^4+r^2 s^7-r s^4+r s^2-1\right)}{s^5} \\
 & \\[-0.6em]
\hline
 & \\[-0.6em]
\partial R_{11} / \partial z & \frac{r^7 s^5-r^4 s^6-r^4 s^5-r^4 s^4-r^4 s^3+r^3 s^4+r^2 s^7+r^2 s^4-r s^5-r s^2-1}{r s^7} \\
 & \\[-0.6em]
\hline
 & \\[-0.6em]
\partial R_{12} / \partial z & -\frac{r^9 s^6-r^7 s^7-r^7 s^4-r^6 s^6-r^6 s^4+r^5 s^5+r^5 s^3+r^5+r s^2-1}{r s^2} \\[-0.6em]
 & \\
\hline
\end{tabular}
\caption{$\partial R_i / \partial z$}\label{tb:ddz}
\end{centering}
\end{table}

\vspace*{\stretch{1}}
\end{center}

\clearpage

\pagebreak

This directly implies Theorem \ref{thm:exists}.

\begin{pf}[Proof of Theorem \ref{thm:exists}]
Let $S_1$ be the Cartwright--Steger surface and $\Gam$ its fundamental group. Then $\Gam^{ab} \cong \bbZ^2$. Let $\al : \Gam \to G$ be a homomorphism onto a finite abelian group of order $n$ and $S_n$ the associated \'etale cover of $S_1$. Theorem \ref{thm:alexander} implies that $b_1(S_n) = b_1(S_1) = 2$. Thus $S_n$ has irregularity $q = 1$. Then $\chi(\calO_{S_1}) = 1$ and $\chi$ is multiplicative in covers, so we have
\[
\chi(\calO_{S_n}) = p_g(S_n) = n.
\]
Then $c_2 = 3 \chi$ for any smooth closed ball quotient, which completes the proof.
\end{pf}

To prove Theorem \ref{thm:moduli}, we must first count the number $a_n(\bbZ^2)$ of subgroups of index $n$ in $\bbZ^2$. In the notation of the proof of Theorem \ref{thm:exists}, this determines the number of distinct $\al : \Gam \to G$ with $G$ a finite abelian group of order $n$. This is well-known to equal $\sig(n)$, where
\[
\sig(n) = \sum_{d\ \mid\ n} d
\]
is the divisor sum function. See \cite[p.\ 308]{LubotzkySegal}. Applying the obvious lower bound $\sig(n) \ge n + 1$, we see that the fundamental group $\Gam$ of the Cartwright--Steger surface has at least $n+1$ normal subgroups of index $n$ with abelian quotient.

To count these surfaces up to homeomorphism, by Mostow--Siu Rigidity \cite{Siu} we must count these subgroups of $\Gam$ up to conjugacy in $\Isom(\bbB^2)$, as opposed to conjugacy in $\Gam$ itself. We now do this to prove Theorem \ref{thm:moduli}.

\begin{pf}[Proof of Theorem \ref{thm:moduli}]
Let $\Gam$ be the fundamental group of the Cartwright--Steger surface and $A(\Gam)$ be the set of equivalence classes of homomorphisms of $\Gam$ onto finite abelian groups, where two homomorphisms are equivalent if they have the same kernel. Given $\al \in A(\Gam)$, let $\Gam_\al$ be the kernel of $\al$. We must show that there is a universal constant $d$ such that $\Gam_\al$ is conjugate in $\Isom(\bbB^2)$ to at most $d$ other $\Gam_\beta$ for $\beta \in A(\Gam)$. Then we can take $c = 1/d$ in the statement of the theorem.

Recall that $\Gam$ is arithmetic. In fact, $\Gam$ is a \emph{congruence subgroup}; see the remark on p.\ 90 of \cite{Stover}. Specifically, $\Gam$ is contained in the arithmetic lattice $\conj{\Gam} = \PSU(2,1; \bbZ[\zeta_{12}])$ in $\PU(2,1)$, where $\zeta_{12}$ is a primitive $12^{th}$ root of unity. The \emph{principal congruence subgroups} of $\conj{\Gam}$ are the finite groups $G(\calI)$ given by taking the image of $\conj{\Gam}$ in the finite group $\PSL_3(\bbZ[\zeta_{12}] / \calI)$ for $\calI$ an ideal of $\bbZ[\zeta_{12}]$. The kernel of this homomorphism will be denoted by $\conj{\Gam}(\calI)$. We note that the groups $G(\calI)$ are perfect groups \cite[\S 6.1]{LubotzkySegal}. (Notice that the groups there are absolutely almost simple and simply connected, whereas ours are adjoint, but this means that our groups are quotients of perfect groups, hence are also perfect.)

The strong approximation theorem \cite[Thm.\ 16.4.2]{LubotzkySegal} implies that $\Gam$ maps onto $G(\calI)$ for all but finitely many ideals $\calI$. Fix one such $\calI$. We claim that none of the subgroups $\Gam_\al$ can contain $\conj{\Gam}(\calI)$. Indeed, this would imply that $\Gam_\al$ would map onto a proper normal subgroup of the perfect group $G(\calI)$ with abelian quotient, which is absurd. For any of the remaining ideals $\calJ$ for which $\Gam$ does not map onto $G(\calJ)$, we see that some $\Gam_\al$ map contain $\conj{\Gam}(\calJ)$, but $\Gam_\al$ then must have bounded index in $\Gam$, hence it follows that only finitely many of the $\Gam_\al$ can be congruence subgroups.

Every arithmetic lattice is contained in finitely many \emph{maximal} arithmetic lattices, and maximal arithmetic lattices are congruence subgroups \cite[Prop.\ 1.4(iv)]{Borel--Prasad}. Since congruence subgroups are closed under intersection, any arithmetic lattice has a well-defined \emph{congruence closure} $C(\Lam)$, the intersection of all the congruence subgroups that contain $\Lam$. Since $\Gam$ is a congruence subgroup, we see that $C(\Gam_\al) \subset \Gam$ for every $\al \in A(\Gam)$. In particular, $C(\Gam_\al) = \Gam_\beta$ for some $\beta \in A(\Gam)$. We showed above that only finitely many $\Gam_\beta$ can be a congruence subgroup, hence
\[
C(A(\Gam)) = \{C(\Gam_\al) : \al \in A(\Gam)\} = \{\Gam_{\al_i} : 1 \le i \le r\}
\]
for some finite subset $\{\al_i\}$ of $A(\Gam)$.

We now consider the set
\begin{align*}
&\calN = \\
&\left\{g \in \mathrm{Comm}(\Gam) : g C(\Gam_\al) g^{-1} = C(\Gam_\beta)\ \textrm{for some}\ C(\Gam_\al), C(\Gam_\beta) \in C(A(\Gam))\right\}.
\end{align*}
Note that $\Gam \subseteq \calN$ since $\Gam$ normalizes each $C(\Gam_\al)$. To prove the theorem, it suffices to prove that the set $\calN / \Gam$ is finite. To see that this does suffice to prove the theorem, suppose that $\calN / \Gam$ has representatives $g_1, \dots, g_d$. If $g \Gam_\al g^{-1} = \Gam_\beta$ for some $\al, \beta \in A(\Gam)$, then $g \in \calN$ by the above. Then $g = g_i \gam$ for some representative $g_i$ and some $\gam \in \Gam$, so
\[
\Gam_\beta = (g_i \gam) \Gam_\al (g_i \gam)^{-1} = g_i \Gam_\al g_i^{-1},
\]
and thus $\Gam_\al$ is conjugate to at most $d$ groups $\Gam_\beta$ for $\beta \in A(\Gam)$. The theorem follows immediately.

Since we already proved that $C(A(\Gam))$ is finite, to prove that $\calN / \Gam$ is finite it suffices to show that
\[
\left\{g \in \mathrm{Comm}(\Gam) : g C(\Gam_\al) g^{-1} = C(\Gam_\beta)\right\}
\]
is finite modulo $\Gam$ for any fixed $\al, \beta \in A(\Gam)$. Indeed, if $g$ conjugates $\Gam_\al$ to $\Gam_\beta$, then it also conjugates $C(\Gam_\al)$ to $C(\Gam_\beta)$. However, every $C(\Gam_\al)$ is of the form $\Gam_\de$ for some $\de \in A(\Gam)$, so it in fact suffices to show that
\[
\calN_{\al, \beta} = \left\{g \in \mathrm{Comm}(\Gam) : g \Gam_\al g^{-1} = \Gam_\beta \right\}
\]
is finite modulo $\Gam$ for any fixed $\al, \beta \in A(\Gam)$.

If $g, h \in \calN_{\al, \beta}$, then $g^{-1} h$ normalizes $\Gam_\al$ in $\Isom(X)$. It is well-known that the normalizer $N(\Gam_\al)$ of $\Gam_\al$ is a lattice in $\Isom(X)$. Moreover, since $\Gam_\al$ is normal in $\Gam$, we see that $\Gam$ is a finite-index subgroup of $N(\Gam_\al)$. Let $n_1, \dots, n_r$ be representatives for $N(\Gam_\al) / \Gam$. Then there is some $n_i$ and a $\gam \in \Gam$ such that $g^{-1} h = n_i \gam$. In particular, $h$ is equal to $g n_i$ modulo the right-action of $\Gam$, and it follows after fixing $g$ and letting $h$ vary over $\calN_{\al, \beta}$ that $\calN_{\al, \beta} / \Gam$ is finite. This completes the proof of the theorem.
\end{pf}

\begin{rem}
The reader may be a bit surprised that the proof of Theorem \ref{thm:moduli} is so involved. However, there is some good reason for the complexity of the argument. One can use the fact that the commensurator $\mathrm{Comm}(\Gam)$ is analytically dense in a finite index subgroup of $\Isom(X)$ (since $\Gam$ is arithmetic), to find large and interesting collections of subgroups of $\Gam$ that are conjugate in $\Isom(X)$ but not in $\Gam$. This leads to the well-studied notion of \emph{hidden symmetries}, and the arithmetic manifolds are precisely those with infinitely many hidden symmetries. See \cite{FS} for more on this. That we can exert so much control on the number of hidden symmetries among our coverings is a consequence of the fact that all our coverings are abelian.
\end{rem}

\begin{rem}
More delicate counting results for arithmetic lattices allow one to also study the case $q=0$. We expect the following to be true. For $n \ge 0$, let $\calZ_n$ be the moduli space of minimal smooth projective surfaces of general type with $q=0$, $p_g = n$, and $c_2 = 3 p_g + 3$. Then $\calZ_n$ is a finite set of points and there is an infinite sequence $\{n_j\}$ such that $\calZ_{n_j}$ is nonempty. Moreover, there is a universal constant $c > 0$ such that for any $k \in \bbZ$ there exists an infinite number of $n$ for which
\[
k \le \#\calZ_n < e^{c n}
\]
for all $n \in \bbZ$.
\end{rem}

\begin{rem}
We close with a final remark on our presentation for the fundamental group of the Cartwright--Steger surface. While the complete details of its construction are unpublished, one can confirm its existence independently of \cite{CS}. As is well-known, this surface is a finite index subgroup of a Deligne--Mostow lattice (e.g., see \cite[p.\ 90]{Stover}), and one can use a presentation for the Deligne--Mostow lattice and Magma \cite{Magma} to find an independent presentation for the fundamental group of the Cartwright--Steger surface. Magma also immediately checks that this presentation is equivalent to the one given by \cite{CS}.
\end{rem}

\bibliography{CScovers}
\end{document}